\documentclass{anstrans}
\title{Developing an Analytical Fixed Source Solver for the 1D Multigroup $S_N$ Equations}
\author{Jilang Miao,$^{*1}$ and Miaomiao Jin$^{*}$}

\institute{
$^{*}$Department of Nuclear Engineering, The Pennsylvania State University, University Park, 16802 PA, USA
}
\footnotetext[1]{\footnotesize Corresponding author: Jilang Miao (jlmiao@psu.edu)}
\usepackage{graphicx} % allows inclusion of graphics
\usepackage{subfigure}

\usepackage{booktabs} % nice rules (thick lines) for tables
\usepackage{microtype} % improves typography for PDF
\usepackage{upgreek}
\usepackage{amsmath}
\usepackage{bbm}
\usepackage{algorithm}
\usepackage{algpseudocode}

\newcommand{\SN}{S$_N$~}
\expandafter\def\expandafter\normalsize\expandafter{%
    \normalsize%
    \setlength\abovedisplayskip{1pt}%
    \setlength\belowdisplayskip{1pt}%
    \setlength\abovedisplayshortskip{1pt}%
    \setlength\belowdisplayshortskip{1pt}%
}

\begin{document}
%%%%%%%%%%%%%%%%%%%%%%%%%%%%%%%%%%%%%%%%%%%%%%%%%%%%%%%%%%%%%%%%%%%%%%%%%%%%%%%%
\section{Introduction}
The discrete ordinates method, commonly known as the \SN method, %,Lar2008
is a discretization of the differential form of the particle transport equation~\cite{Car1953}. 
The angular flux is then solved at specific angles, i.e., this method relies on a conceptually straightforward evaluation of the transport equation at a limited number of discrete angular directions, or ordinates. 
Additionally, quadrature relationships are employed to replace integrals over angle, simplifying the integration with summations over these discrete ordinates~\cite{hebert2009applied}. %,stacey2018nuclear

Earlier work developed an accurate eigenvalue solver for multigroup \SN equations in slab geometry. In this method, analytical \SN solution was obtained for each homogeneous subregion by expansion on eigensystem which is determined by neutron cross sections in the material~\cite{analytical2A2G,AMGSN2023}. 
The expansion coefficients were solved from a linear system incorporating the continuity condition at the interfaces and boundary condition of the angular fluxes. 
The eigenvalues were found by searching the root of the determinant of the boundary condition matrix. 

In this work, we intend to extend the analytical multigroup \SN method to solve 1D fixed source problems.
The 1D fixed source \SN solution can be  applied as the axial solver in the 2D-1D schemes to simulate 3D transport where the radial traverse leakage is viewed as external source for the 1D problem.
Similarly, the fixed source solution can be used to develop 3D nodal \SN methods.
In addition, the 1D fixed source solution are commonly used in the iterative methods such as power iteration to solve 1D eigenvalue problems, where determinant root finder becomes intractable for large matrices.

To approach this, we first derive the solution for the 1D fixed source problem based on the analytical multigroup \SN method. 
Then the solver is applied in the power iteration for an 1D eigenvalue problem.
%The test case is created by generating multigroup cross sections for a pincell with OpenMC~\cite{romano2013openmc,boyd2019multigroup}.
%The fixed source solver will then find the eigenvalue and eigenfunction of the slab problem.
For performance study, we also use the traditional sweeping based \SN algorithm to solve the problem and
demonstrate the fixed source solver is highly accurate and efficient.

%%%%%%%%%%%%%%%%%%%%%%%%%%%%%%%%%%%%%%%%%%%%%%%%%%%%%%%%%%%%%%%%%%%%%%%%%%%%%%%%
\section{Theory}
\subsection{S$_N$ equation in a homogeneous slab}

For a given number of energy groups, denoted as $g=1,...,G$, 
and a quadrature set $\left . \{\mu_n,\omega_n\} \right | _{n=1,...,N}$, 
the transport equation for the angular flux $\psi_{g,n}$ is expressed in Eq.~\ref{eq::sn1d}. 
\begin{equation}
\begin{aligned}
 \mu_n \frac{\partial}{\partial x} \psi_{g,n}(x)+  \Sigma_{t, g} \psi_{g,n}(x)  = 
   \sum_{n^{\prime}, g^{\prime}} \omega_{n \prime} \Sigma_{s, g^{\prime} n^{\prime} \rightarrow g n} \psi_{g^{\prime}, n^{\prime}}(x) \\
+
%\frac{1}{k_{eff}}
\sum_{n^{\prime}, g^{\prime}} \omega_{n^{\prime}}  \nu \Sigma_{f, g^{\prime} n^{\prime} \rightarrow g n} \psi_{g^{\prime}, n^{\prime}}(x) 
+Q_{g, n}(x)
\label{eq::sn1d}
\end{aligned}
\end{equation}

The angular flux $\psi_{g,n}$ can be compactly aggregated in a vector $\Psi(x)$ of length $N G$. 
This vector consists of $G$ blocks, each having a length of $N$. For a specific block $g$ ($g=1,...,G$), 
it corresponds to the angular fluxes $\left . \psi_{g,n} \right|_{n=1,\cdots,N}$. 
Consequently, we can denote $\psi_{g,n}(x)$ as $\Psi_{g N + n}(x)$.
Similarly, the source $Q$ can be represented in a $N G$ vector $Q(x)$, such that $Q_{g,n}(x)=Q_{g N + n}(x)$.

With the same convention as in \cite{AMGSN2023} to organize the cross-sections and quadrature sets into matrices,
Eq.~\ref{eq::sn1d} can be written in matrix form as in Eq.~\ref{eq::sn1d matrix}.
\begin{equation}
    \partial_x \Psi(x) = A \Psi(x) + \Theta(x)
    \label{eq::sn1d matrix}
\end{equation}
where 
\begin{equation}
     \Theta_{g N + n}(x) = \frac{1}{\mu_n} Q_{g N + n}(x) 
\label{eq::ThetaQ}
\end{equation}
The solution to Eq.~\ref{eq::sn1d matrix} through the process of block- diagonalization \cite{strang2006linear} of the matrix $A$ .
In general, a real matrix $A$ is similar to a block-diagonal matrix $B$, i.e., there exists an invertible matrix $P$ such that 
\begin{equation}
A P = P B    .
\label{eq::AP=PB}
\end{equation}

Since the $\partial_x$ operator commutes with a constant matrix ($P^{-1}$), after applying $P^{-1}$ to both sides of
the equation for $\Psi(x)$ (Eq.~\ref{eq::sn1d matrix}), the following can be obtained, 
\begin{equation}
    \partial _x X(x)  = B X(x) +  P^{-1}\Theta(x)
    \label{eq:sX=BX+Q}
\end{equation}
where 
\begin{equation}
    X(x) := P^{-1} \Psi(x)
    \label{eq::X=invP Psi}
\end{equation}

With $B$ being a block-diagonal matrix, we can derive the solution for $X(x)$ in Eq~\ref{eq:sX=BX+Q} 
(derivation details will be provided in the full article), 
\begin{equation}
    X(x) = \Gamma(x) \left (  \alpha + 
    \int_{x_0}^x d\xi \Gamma(-\xi) P^{-1} \Theta(\xi)
    \right)
    \label{eq::X=G a}
\end{equation}
where $\alpha$ is an undetermined vector to be solved. 
The details on constructing the block-diagonal 
 matrices $B$ and $\Gamma$ from the eigensystem of $A$ can be found in \cite{AMGSN2023}. Particularly, $\Gamma$ contains exponential and trigonometric functions, and hence, 
the integral term in Eq.~\ref{eq::X=G a} can be analytically computed 
for a broad spectrum of functions for the source $\Theta$, such as polynomials, exponentials, and trigonometric functions.

Specifically, if the external source $Q$ is piece-wise constant over a mesh comprising $M$ regions, i.e.,
\begin{equation}
    Q(x) = Q_{m}\mathbbm{1}_{x \in[x_{m-1},x_m)}, m = 1,\dots, M
    \label{eq::piecewise Q}
\end{equation}
Substituting the piece-wise constant $Q$ into Eq.~\ref{eq::X=G a}, 
we obtain a more simplified expression for the solution of $X(x)$,
\begingroup
\scriptsize
\begin{equation}
\begin{split}
    X(x) = \Gamma(x)  (  \alpha 
    %& 
    + 
    \sum_{m=1}^{m^*}\left[
        \Gamma(-x_{m-1}) - \Gamma(-x_m)
    \right] B^{-1} P^{-1} \Theta_m 
    %\\
    %&
     + \left[
        \Gamma(-x_{m^*}) - \Gamma(-x)
    \right] B^{-1} P^{-1} \Theta_{m^*}
    )
\end{split}
    \label{eq::X=G a piecewise}
\end{equation}
\endgroup
where $m^*$ denotes the index of the mesh grid containing $x$.

Finally, it is noted that the selection of $x_0$ in Eq.~\ref{eq::X=G a} and Eq.~\ref{eq::X=G a piecewise} is arbitrary, 
as the integral constant can be combined into $\alpha$. 
For convenience, it could be chosen to be the left boundary of the region.

\subsection{S$_N$ solution in a heterogeneous slab}
\label{sec::sn sln}

Consider a heterogeneous slab which can be divided into $R$ individually homogeneous regions, numbered as $1,\cdots,R$ from left to right. 
The position of the separating interfaces are defined as $x_0,\cdots,x_R$.
Hence, we need to determine the $\alpha$ term in Eq.~\ref{eq::X=G a} for each region. To achieve this, we resort to i) boundary conditions (left and right) and ii) continuity requirements for the angular fluxes at region interfaces. 

First, we will formulate the equations corresponding to boundary conditions of incoming fluxes (e.g., vacuum boundaries as zero incoming fluxes).%and reflective boundaries.
The formulation for reflective boundary conditions does not depend on the external source and can be found in~\cite{AMGSN2023}.
Then, we will consolidate the equations for each homogeneous region into one linear system to solve for the solution.
Here, $\alpha_i$ (vector of length $NG$) denotes the coefficients for region $i$.
The coefficients for all regions will be consolidated into vector $\alpha$, with a length of $NGR$.
Similarly, $P_i$ and $\Gamma_i$ represent the transform matrix $P$ and block-diagonal matrix $\Gamma$ for region $i$, respectively.

\subsubsection{Incoming flux boundary condition}
For incoming flux $\Psi_L$ from the left end, i.e., angular flux with $\mu>0$, the boundary condition can be represented as
\begin{equation}
    \left. \left( P_{1} \Gamma_{1} (x_0) \right)\right| _{\left\{\mu>0\right\}} \alpha_1 = \Psi_L ,
    \label{eq::BC srcL}
\end{equation}
%based on Eq.~\ref{eq::X=G a}.

Similarly, for incoming source $\Psi_R$ from the right end, the boundary condition can be represented as
\begingroup
\scriptsize
\begin{equation}
\begin{split}
  \left. \left( P_{R} \Gamma_{R} (x_R) \right)\right| _{\left\{\mu<0\right\}} \alpha_R  = \Psi_R
  -
  \left. \left( P_{R} \Gamma_{R} (x_R) \right)\right| _{\left\{\mu<0\right\}} 
    \int_{x_{R-1}}^{x_R} d\xi \Gamma_R(-\xi) P_R^{-1} \Theta(\xi)
\label{eq::BC srcR}
\end{split}
\end{equation}
\endgroup

In Eq.~\ref{eq::BC srcL} and Eq.~\ref{eq::BC srcR}, 
$\left. M \right|_{\{\mu>0\}} $ (or $\left. M \right|_{\{\mu<0\}}$) means the operation of extracting specific rows from matrix $M$, where $M$ is a placeholder for the $P\Gamma$ matrices ($NG\times NG$) noted in the equations. 
This selection is based on the following procedure: the discrete angles $\left\{\mu_n\right\}_{n=1,\cdots,N}$ are repeated $G$ times to form a vector of length $NG$, and then the rows corresponding to $\mu>0$ (or $\mu<0$) are selected. 

\iffalse
\subsubsection{Reflective boundary condition}
The reflective boundary on the left end can be written as,  
\begin{equation}
   \left[
   \left.\left( P_{1} \Gamma_{1} (x_0) \right)\right|_{\left\{\mu<0\right\}} -
   \left.\left( P_{1} \Gamma_{1} (x_0) \right)\right|_{\left\{\mu>0\right\}}
   \right]\alpha_1 = 0  
   \label{eq::BC refL}
\end{equation}
and the reflective boundary on the right end can be written as,  
\begin{equation}
   \left[
   \left.\left( P_{R} \Gamma_{R} (x_R) \right)\right|_{\left\{\mu<0\right\}} -
   \left.\left( P_{R} \Gamma_{R} (x_R) \right)\right|_{\left\{\mu>0\right\}}
   \right]\alpha_R = 0  
   \label{eq::BC refR}
\end{equation}
In Eq.~\ref{eq::BC refL} and Eq.~\ref{eq::BC refR}, 
when performing the row selections for $\left. M \right|_{\{\mu>0\}}$ and $\left. M \right|_{\{\mu<0\}}$, 
it is important to ensure that for the same row, the corresponding $\mu$ values have the same absolute values. 
Such requirement is automatically satisfied when the discrete angles $\left\{\mu_n\right\}_{n=1,\cdots,N}$ are pre-sorted in monotonically increasing/decreasing absolute value. 
\fi 

\subsubsection{Angular flux continuity condition}
At region interfaces, all angular fluxes are continuous.
Hence, the condition for the interface between region $i$ and $i+1$ at $x_i$ can be written as 
\begingroup
\footnotesize
\begin{equation}
    P_{i} \Gamma_{i} (x_i) \alpha_i  -
     P_{i+1} \Gamma_{i+1} (x_i)   \alpha_{i+1}
    = - P_{i} \Gamma_{i} (x_i) \int_{x_{i-1}}^{x_i} d\xi \Gamma_i(-\xi) P_i^{-1} \Theta(\xi)
   \label{eq::BC cont}
\end{equation}
\endgroup

\subsubsection{Solution of the coefficients}
With the boundary conditions specified in Eqs.~\ref{eq::BC srcL}-
%\ref{eq::BC refR} 
\ref{eq::BC srcR} 
for the two ends, each end yields $NG/2$ equations, leading to a total of $NG$ equations.
At the $R-1$ interior interfaces, the continuity requirement leads to $(R-1)NG$ equations. In sum, there are $NG\times R$ equations that can be combined to solve for the $NG\times R$ coefficients in $\left. \left\{\alpha_i\right\} \right|_{i=1,...,R}$. For ease of notation, we have following definitions,
\begin{align}
P\Gamma_1^+&\equiv \left. \left( P_1 \Gamma_1(x_0) \right)\right|_{\left\{\mu>0\right\}}  \\
%P\Gamma_1^-&\equiv \left. \left( P_1 \Gamma_1(x_0) \right)\right|_{\left\{\mu<0\right\}}  \\
%P\Gamma_R^+&\equiv \left. \left( P_R \Gamma_R(x_R) \right)\right|_{\left\{\mu>0\right\}}   \\
P\Gamma_R^-&\equiv \left. \left( P_R \Gamma_R(x_R) \right)\right|_{\left\{\mu<0\right\}}      
\end{align}

As a specific example, considering the boundary condition where both ends have incoming sources, 
the coefficients $\left. \left\{\alpha_i\right\} \right|_{i=1,...,R}$ can be determined through Eq.~\ref{eq::BCmat VV}, which consists of a linear system with dimension $NGR\times NGR$. The boundary conditions at both ends are placed in the first $NG$ rows of this arrangement.
The interface conditions are placed in the remaining $(R-1)NG$ rows, 
with each interface contributing $NG$ rows. 
The matrix on the left-hand side is structured as a block matrix of size $R\times R$, 
where each block represents a matrix of dimension $NG\times NG$. 
The vector on the right-hand side is presented as a $R\times 1$ vector, 
with each entry itself is a vector with length $NG$. 
\begingroup
\footnotesize
\begin{equation}
\begin{split}
\left[\begin{array}{cccc}
\left[\begin{array}{c} 
P\Gamma_1^+ 
%define this as P\Gamma_1^+
\\ \mathbf{0} 
\end{array}
\right]
& \mathbf{0} & \ldots  & 
\left[\begin{array}{c} 
\mathbf{0} \\ 
P\Gamma_R^-
\end{array}
\right]
\\
P_1 \Gamma_1(x_1) & - P_2 \Gamma_2(x_1) & \ldots &  \mathbf{0} \\
\ldots & \ldots & \ldots &   \ldots \\
\mathbf{0} &  \ldots & P_{R-1}\Gamma_{R-1}(x_R) & -P_R\Gamma_R(x_R)
\end{array}\right] \alpha \\ = 
\left[\begin{array}{c}
\left[
\begin{array}{c} 
\Psi_L \\ \Psi_R
    -  P\Gamma_R^- 
    \int_{x_{R-1}}^{x_R} d\xi \Gamma_R(-\xi) P_R^{-1} \Theta(\xi)
\end{array}
\right]  
 \\
- P_{1} \Gamma_{1} (x_1) \int_{x_{0}}^{x_1} d\xi \Gamma_1(-\xi) P_1^{-1} \Theta(\xi) \\
\ldots \\
- P_{R-1} \Gamma_{R-1} (x_{R-1}) \int_{x_{R-2}}^{x_{R-1}} d\xi \Gamma_{R-1}(-\xi) P_{R-1}^{-1} \Theta(\xi)
\end{array}\right]
\label{eq::BCmat VV}
\end{split}
\end{equation}
\endgroup

The linear systems for other boundary conditions can be constructed in a similar way and are skipped in this summary.
To acquire the matrices $P$ and $\Gamma$ as used in Eq.~\ref{eq::BCmat VV}, 
it is necessary to determine the eigensystem of matrix $A$ (in Eq.~\ref{eq::sn1d matrix}) for each region.
Since $A$ only depends on cross-sections of the material, if the $R$ regions span $M$ distinct materials ($M\leq R$), it is only necessary to find $M$ such eigensystems.
Complexity to construct the matrices ($P$ and $\Gamma$ in Eq~\ref{eq::BCmat VV}) is thus on the order of $M\times \mathcal{O}( (NG)^3 )$, where $\mathcal{O}( (NG)^3 )$ is from sovling the eigensystem. The complexity to solve the $NGR\times NGR$ linear system (Eq~\ref{eq::BCmat VV}) is $\mathcal{O}( (NGR)^3 )$ based on matrix inversion. 
%For an eigenvalue problem, the right side of
%Eq.~\ref{eq::BCmat VV} and Eq.~\ref{eq::BCmat RV} 
%becomes zero; hence $k_{eff}$ can be found by solving the determinant of the matrix on left hand equal to $0$. 
%The roots of the determinant yield all eigenmodes of the problem. 
%The root-finding algorithm is expected to identify the largest eigenvalue, $k_{eff}$, corresponding to the fundamental mode.
%Given the $k_{eff}$, the coefficients $\alpha$ can be obtained by SVD decomposition of the matrix, 
%which has a complexity of $\mathcal{O}\left( (NGR)^3 \right)$.
%It should be noted that, in the process of searching $k_{eff}$, the $M$ eigensystems need to be updated for each $k_{eff}$.

\subsection{Application of the fixed source solution}
In this section, we apply the fixed source solver to eigenvalue problems.
Especially, if power iteration is used to find the fundamental mode,
each iteration step corresponds to a fixed source problem.
In iteration  $n$,
fission term in Eq.~\ref{eq::sn1d} can be treated as the external source,
\begin{equation}
Q^{(n)}(x) = \left( \frac{1}{k^{(n)}_{eff}} - \frac{1}{k_e} \right ) \sum_{n^{\prime}, g^{\prime}} \omega_{n^{\prime}}  \nu \Sigma_{f, g^{\prime} n^{\prime} \rightarrow g n} \psi^{(n)}_{g^{\prime}, n^{\prime}}(x) 
\label{eq::Qn}
\end{equation}
%Hence, the matrix $A$ (in Eq.~\ref{eq::sn1d}) does not have contributions from fission cross-sections.
For acceleration purposes in the power iteration,
the solver here has the flexibility
to allow Wielandt's shift in $k_e$ ~\cite{brown2007wielandt}. 
Notably, the block-diagonalization of matrix $A$ in this work can efficiently treat the complex eigenvalues of $A$ resulting from Wielandt's shift.

With the source term in Eq.~\ref{eq::Qn} reasonably assumed using piece-wise constant functions on a fine mesh with size $M$ ($M$>>$R$), Eq.~\ref{eq::X=G a piecewise} can be used to calculate the integral required to solve $\alpha$ vector. The corresponding algorithm is summarized in Algorithm~\ref{alg::power}.
\begingroup
\footnotesize
\begin{algorithm}
  \caption{Eigenvalue power iteration with fixed source analytical multigroup \SN solver}
  \label{alg::power}
\begin{algorithmic}

\For{each distinct material}
\State construct matrix $A$
\State find block-diagonalization matrices $P$ and $B$ 
\EndFor
  \State initialize piecewise constant fission source $Q^{(0)}$ % \Comment{vector length is $MNG$}

\While{ error metric above threshold }
    \State solve coefficients $\alpha^{(n)}$ (Eq.~\ref{eq::BCmat VV}) %~\ref{eq::BCmat RV}~\ref{eq::BCmat RR})
    \State evaluate $\Psi^{(n)}$ on the source mesh centers (Eq.~\ref{eq::X=G a piecewise} )
    \State calculate  $Q^{(n)}$ from  $\Psi^{(n)}$
    \State update $k^{(n)}_{eff}$
    \State calculate error metric such as norm of $Q^{(n)}-Q^{(n-1)}$
%    \State $n \gets n + 1$
\EndWhile
\end{algorithmic}
\end{algorithm}
\endgroup

Note that, here a fine mesh is used to describe the source term based on piece-wise constant functions, while the whole system is still described by the $R$ homogeneous regions on a coarse mesh (the linear system is of size  $NGR\times NGR$). For the case where i) one energy group is assumed, ii) angular fluxes are solved on the same fine mesh as source $Q$, and iii) there is no Wielandt's shift in $Q$, 
Algorithm~\ref{alg::power} is reduced to the earlier work in~\cite{wang2017new}.

%%%%%%%%%%%%%%%%%%%%%%%%%%%%%%%%%%%%%%%%%%%%%%%%%%%%%%%%%%%%%%%%%%%%%%%%%%%%%%%%
\section{Results}

As a test case, we study a 35 cm slab with 3 regions.
The reactor core is located within [-15 cm, 15 cm].
The reflector is within [-17.5 cm, -15 cm] and [15 cm, 17.5cm].
The system has vacuum boundary condition on both ends.
Two-group cross-sections (in the unit of cm$^{-1}$) for the core and reflector materials are shown in Table~\ref{tab::xs}, which are generated with OpenMC~\cite{romano2013openmc,boyd2019multigroup} for a typical fuel pincell. 
\begingroup
\scriptsize
\begin{table}[htb]
  \centering
  \caption{Cross-section parameters.}
  \begin{tabular}{llr}\toprule
      &  Core      & Reflector
\\ \midrule
$\Sigma_{t,1}$              & 6.8294e-01 & 8.9176e-01 \\
$\Sigma_{t,2}$              & 2.0658e+00 & 3.0361e+00 \\
$\Sigma_{s,1\rightarrow 1}$ & 6.4870e-01 & 8.4530e-01 \\
$\Sigma_{s,1\rightarrow 2}$ & 2.5869e-02 & 4.6078e-02 \\
$\Sigma_{s,2\rightarrow 1}$ & 4.2114e-04 & 2.8498e-04 \\
$\Sigma_{s,2\rightarrow 2}$ & 1.9696e+00 & 3.0181e+00 \\
$\nu\Sigma_{f,1}$           & 6.0427e-03 & 0.0000e+00 \\
$\nu\Sigma_{f,2}$           & 1.5343e-01 & 0.0000e+00 \\
$\chi_1$                    & 1.0000e+00 & 0.0000e+00 \\
$\chi_2$                    & 0.0000e+00 & 0.0000e+00 \\ 
\bottomrule
\end{tabular}
  \label{tab::xs}
\end{table}
\endgroup

A reference solution is generated using OpenMC~\cite{romano2013openmc} multigroup mode with
the same geometric configuration, boundary conditions and cross-sections. 
The simulation tracks $10^6$ neutrons per generation.
The neutrons are simulated for $200$ inactive generations and tallies are collected for the next $800$ active generations to compute  
scalar fluxes, angular fluxes and $k_{eff}$.
The fluxes are tallied on a spatially uniform mesh of size $700$ for each energy group.
In addition, the angular fluxes are tallied over a specific polar angle range corresponding to the \SN quadrature set.

\subsection{Accuracy of the eigenvalue problem}
With Gauss-Legendre quadrature sets,
Algorithm~\ref{alg::power} is used to run the power iteration for $S_2$, $S_4$, $S_8$ and $S_{16}$.
The initial guess of the source term is isotropic and varies according to $Q_{g,n}(x) \propto |x|$ .
The iteration is terminated when the $L^2$ norm of scalar flux ($\phi$) change between two consecutive generations is below $10^{-6}$.
\begin{equation}
||\phi^{(n)}-\phi^{(n-1)}||_2 < 10 ^{-6}
\label{eq::exerr}
\end{equation}
%For the two-group test problem on $700$ mesh, $\phi^{(n)}$ are vectors of length $1400$.
We note that for all orders, the solution converges after around $25$ iterations.
To compare with Monte Carlo (MC) reference, 
the fluxes from \SN are normalized such that the sum of the integral of the scalar fluxes over all groups is 1. Table~\ref{tab::keff} shows the $k_{eff}$ from OpenMC and the different orders of the analytical \SN solvers.
It clearly shows how higher order solution approaches the MC reference.

\begingroup
\scriptsize
\begin{table}[htb]
  \centering
  \caption{Computed $k_{eff}$ compared with MC reference.}
  \begin{tabular}{llr}\toprule
Method      & $k_{eff}$      & $k_{eff}$ - $k_{eff,MC}$ (pcm) 
\\ \midrule
MC reference & 1.24953 $\pm$ 0.00002 &  \\
\midrule
Analytical $S_2$        & 1.24737 & -216 \\
Analytical $S_4$        & 1.24936 & -17 \\
Analytical $S_8$        & 1.24949 & -4 \\
Analytical $S_{16}$     & 1.24952 & -1 \\
\midrule
Sweeping   $S_2$        & 1.24288 & -665 \\
Sweeping   $S_4$        & 1.24536 & -417 \\
Sweeping   $S_8$        & 1.24562 & -391 \\
Sweeping   $S_{16}$     & 1.24569 & -384 \\
    \bottomrule
\end{tabular}
  \label{tab::keff}
\end{table}
\endgroup

Next, we proceed to compare the scalar fluxes.
%Fig.~\ref{fig::phi} 
Fig. 1(a--h) present the comparison,
including the results from $S_2$, $S_4$, $S_8$ and $S_{16}$.
%In each subfigure, 
In Fig. 1(a \& e), the scalar fluxes from $S_{16}$ and MC are compared  for fast and thermal group, respectively.
The upper plots show the accurate match of the scalar fluxes,
and the bottom plots indicate the point-wise relative error between S$_{16}$ and MC reference is around $0.75\%$ and $0.1\%$ for fast group and thermal group, respectively. 
%Fig.~\ref{fig::phiS2} demonstrates that $S_2$ can already yield an approximately correct shape.
The point-wise relative error decreases from around $10\%$ in $S_2$ to around $0.1\%$ in $S_{16}$. 
Hence, with increasing orders, a drastic improvement in performance is achieved.
Similar conclusions for angular fluxes ($\omega_n \psi_{g,n}$) can be made. 
As shown in 
%Fig.~\ref{fig::psi},
Fig. 1(i--p), 
the angular fluxes from \SN match MC results very well 
and the point-wise relative error 
decreases from around $30\%$ in $S_2$ to around $0.5\%$ in $S_{16}$.

\iffalse
\begin{figure}[htbp]
\centering
\subfigure
{\includegraphics[width=0.2\textwidth]{figures/fixedSourceIntro/pctPhiVacFine.fmfixsrcMCCompL0S2R14.1_12_1.tot700.ErrEx.pdf}
\label{fig::phiS2}} 
\subfigure
{\includegraphics[width=0.2\textwidth]{figures/fixedSourceIntro/pctPhiVacFine.fmfixsrcMCCompL0S4R14.1_12_1.tot700.ErrEx.pdf}
\label{fig::phiS4}}  
\subfigure
{\includegraphics[width=0.2\textwidth]{figures/fixedSourceIntro/pctPhiVacFine.fmfixsrcMCCompL0S8R14.1_12_1.tot700.ErrEx.pdf}
\label{fig::phiS8}}  
{\includegraphics[width=0.2\textwidth]{figures/fixedSourceIntro/pctPhiVacFine.fmfixsrcMCCompL0S16R28.2_24_2.tot700.ErrEx.pdf}
\label{fig::phiS16}}  
\caption{Scalar fluxes from $S_N$ vs. MC. }
\label{fig::phi}
\end{figure}
\fi

\iffalse
\begin{figure}[htbp]
\centering
\subfigure
{\includegraphics[width=0.2\textwidth]{figures/fixedSourceIntro/pctPsiVacFine.fmfixsrcMCCompL0S2R14.1_12_1.tot700.ErrEx.pdf}
\label{fig::psiS2}} 
\subfigure
{\includegraphics[width=0.2\textwidth]{figures/fixedSourceIntro/pctPsiVacFine.fmfixsrcMCCompL0S4R14.1_12_1.tot700.ErrEx.pdf}
\label{fig::psiS4}}  
\subfigure
{\includegraphics[width=0.2\textwidth]{figures/fixedSourceIntro/pctPsiVacFine.fmfixsrcMCCompL0S8R14.1_12_1.tot700.ErrEx.pdf}
\label{fig::psiS8}}  
{\includegraphics[width=0.2\textwidth]{figures/fixedSourceIntro/pctPsiVacFine.fmfixsrcMCCompL0S16R28.2_24_2.tot700.ErrEx.pdf}
\label{fig::psiS16}}  

\caption{Angular fluxes from $S_N$ ($w_n \psi_{g,n}$) vs. MC. }
\label{fig::psi}
\end{figure}
\fi

\begin{figure*}
\centering
  \includegraphics[width=1.0\textwidth]{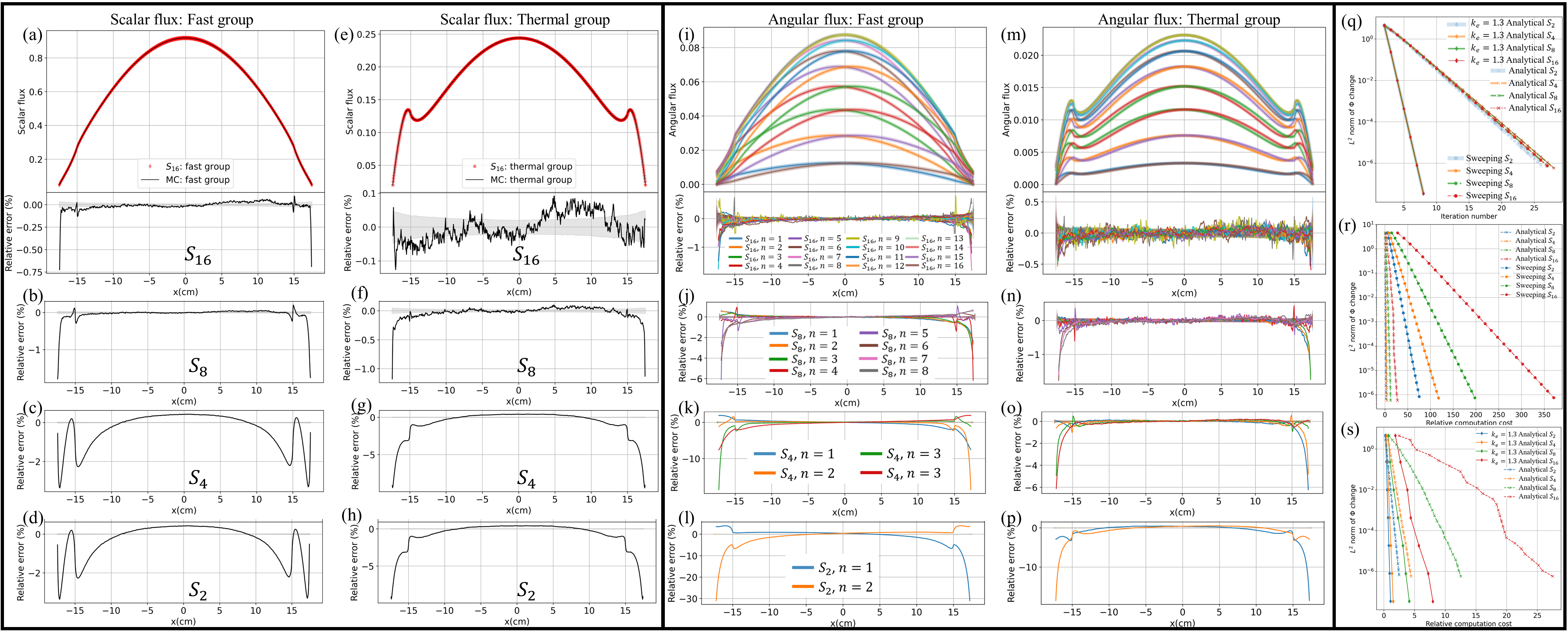}
  \caption{ (a--h) Scalar flux from \SN compared with MC.
  (i--p) Angular flux ($\omega_n \psi_n(x)$) from \SN compared with MC. 
  (a),(e),(i),(m) show both flux value and relative error (\%) from MC. 
  The uncertainty of each tally $T$ from MC is shown with the shading area between $\pm 100\times \frac{\sigma_{\bar{T}}}{\bar{T}}$.
  (q--s) Convergence rate compared with sweeping method. 
  (q) $L^2$ norm of scalar flux change as function of iteration number.
  (r--s) $L^2$ norm of scalar flux change as function of computation time.
  }
\end{figure*}

Further, we compare the accuracy of this solver with the traditional sweeping based \SN method.
The power iteration in the sweeping method is terminated by the same criteria as in Eq.~\ref{eq::exerr}.
The sweeping method requires another layer of iteration for the fixed source problem,
where the inner layer iteration is terminated at half the threshold of outer layer power iteration. 
The comparison of $k_{eff}$ from the two methods are given in Table~\ref{tab::keff}.
The sweeping method indicates significantly larger error than the analytical \SN method.
The reason is that, although both  methods are solving on the same fine mesh of size $700$,
the sweeping method assumes constant source and fluxes in each region, however,
the analytical method only assumes constant source, while the fluxes are analytically represented by eigensystem expansions.

\subsection{Efficiency of the fixed source solver}
Here, we demonstrate the efficiency advantage of the analytical \SN method.
%Fig~\ref{fig::conv iter} 
Fig. 1(q)
plots the $L^2$ norm of scalar flux changes versus number of iterations.
It shows that both the analytical method (without Wielandt's shift) and sweeping method converge at the same rate at all the \SN orders.
They all converge with the same criteria (Eq.~\ref{eq::exerr}) after around $25$ iterations.
This implies that different orders of \SN methods have dominance ratios close to each other despite the $k_{eff}$ differences. 
%Fig~\ref{fig::conv iter} 
Fig. 1(q)
also shows that with the Wielandt's shift $k_e=1.3$, 
the analytical method is significantly accelerated and converges within $10$ iterations. 

We then analyze the computation cost for each iteration. 
%Fig~\ref{fig::conv time} 
Fig. 1(r)
plots the $L^2$ norm of scalar flux change versus time, which is measured in the unit of the average time of solving one iteration in the case of analytical $S_{16}$.
It shows that the analytical method is significantly faster than the sweeping method.
With the same convergence criteria, the analytical method has $31$x,$21$x,$18$x,$22$x speed up for the $S_2$,$S_4$,$S_8$,$S_{16}$ orders, respectively.
Considering that the $S_{16}$ sweeping method has $384$pcm eigenvalue error, which has already been outperformed by the $S_2$ analytical method with $216$pcm eigenvalue error (Table~\ref{tab::keff}),
the analytical method has over $147$x speed up. 
Moreover, with the flexibility of applying the Wielandt's shift in the analytical method, 
%Fig~\ref{fig::conv time w} 
Fig. 1(s)
demonstrates the further improvement in speedup;
Wielandt's factor $k_e=1.3$ largely reduce the number of iterations (from 25 to 10), 
and the treatment of resultant complex eigenvalues in matrix $A$ does not compromise the advantage in computation time. 

\iffalse
\begin{figure}[htbp]
\centering
\subfigure
{\includegraphics[width=0.15\textwidth]{figures/fixedSourceIntro/convergeVsIteration.Cs.0_1_2.png}
\label{fig::conv iter}} 
\subfigure
{\includegraphics[width=0.15\textwidth]{figures/fixedSourceIntro/convergeVsTime.Cs.1_2.png}
\label{fig::conv time}}
\subfigure
{\includegraphics[width=0.15\textwidth]{figures/fixedSourceIntro/convergeVsTime.Cs.0_1.png}
\label{fig::conv time w}}  

\caption{Convergence of analytical and sweeping based \SN methods}
\label{fig::conv}
\end{figure}
\fi

%%%%%%%%%%%%%%%%%%%%%%%%%%%%%%%%%%%%%%%%%%%%%%%%%%%%%%%%%%%%%%%%%%%%%%%%%%%%%%%%
\section{Conclusions}
In this work, we developed the fixed source capability of the analyical multigroup \SN equations in slab geometry.
We demonstrated the application of the fixed source capability in the eigenvalue power iterations.
For the slab problem homogenized from a typical pincell, we observe $216$pcm eigenvalue accuracy for $S_2$ solution
and $1$pcm eigenvalue accuracy in $S_{16}$ solution.
High accuracy was also observed in angular fluxes.
Compared to the sweeping based \SN methods, the analytical method has around $20$x speed up to converge the scalar flux and
around $150$x speedup to reach the same eigenvalue accuracy.

%In the future work, 
%we will extend the 1D solver to 3D neutron transport by schemes like $2D-1D$ coupling and $3D$ nodal methods. 
%In terms of 1D eigenvalue solver, we will develop more advanced iteration method that does not require fine mesh to justify the
%assumption of piecewise constant source. 

%%%%%%%%%%%%%%%%%%%%%%%%%%%%%%%%%%%%%%%%%%%%%%%%%%%%%%%%%%%%%%%%%%%%%%%%%%%%%%%%

%%%%%%%%%%%%%%%%%%%%%%%%%%%%%%%%%%%%%%%%%%%%%%%%%%%%%%%%%%%%%%%%%%%%%%%%%%%%%%%%
\section{Acknowledgments}
This work is supported by the Department of Nuclear Engineering, The Pennsylvania State University.

%%%%%%%%%%%%%%%%%%%%%%%%%%%%%%%%%%%%%%%%%%%%%%%%%%%%%%%%%%%%%%%%%%%%%%%%%%%%%%%%

\footnotesize
\bibliographystyle{ans}
\bibliography{bibliography}
\end{document}